\documentclass[12pt]{article}

\usepackage{amsmath,amssymb,amsfonts,latexsym,theorem}

\numberwithin{equation}{section}
\newtheorem{theorem}{Theorem}[section]

\newtheorem{conjecture}[theorem]{Conjecture}
\newtheorem{lemma}[theorem]{Lemma}
\newtheorem{fact}[theorem]{Fact}

{\theorembodyfont\rmfamily
\newtheorem{definition}[theorem]{Definition}
\newtheorem{remark}[theorem]{Remark}
\newtheorem{example}[theorem]{Example}
\newtheorem{notation}[theorem]{Notation}
}

\def\D{\frak D}
\def\DD{\mathcal D}
\def\S{\frak S}

\def\qed{\hfill $\Box$}

\newenvironment{proof}{\noindent{\itshape Proof.}}{\qed}

\title{\textbf{Restricted Dumont permutations}}

\author{
\textbf{Alexander Burstein}\\
Department of Mathematics\\
Iowa State University\\
Ames, IA 50011-2064 USA}

\date{\today}

\begin{document}

\maketitle

\begin{abstract}
We analyze the structure and enumerate Dumont permutations of the
first and second kinds avoiding certain patterns or sets of
patterns of length 3 and 4. Some cardinalities are given by
Catalan numbers, powers of 2, little Schr\"oder numbers, and other
known or related sequences.
\end{abstract}

\section{Preliminaries} \label{sec:prelim}

\subsection{Patterns}
\label{subsec:patterns}

Let $\sigma\in\S_n$ and $\tau\in\S_k$ be two permutations. We say
that $\sigma$ \emph{contains} $\tau$, or $\tau$ \emph{occurs} in
$\sigma$, if $\sigma$ has a subsequence
$(\sigma(i_1),\dots,\sigma(i_k))$, $1\le i_1<\dots<i_k\le n$,
order-isomorphic to $\tau$. Such a subsequence is called an
\emph{occurrence} (or an \emph{instance}) of $\tau$ in $\sigma$.
In this context, the permutation $\tau$ is called a
\emph{pattern}. We say that $\sigma$ \emph{avoids} $\tau$, or
$\sigma$ is \emph{$\tau$-avoiding}, if $\tau$ does not occur in
$\sigma$.

\begin{notation} \label{not:avoid}
We denote the set of permutations in $\S_n$ avoiding a pattern
$\tau$ by $\S_n(\tau)$. If $T$ is a set of patterns, then we
denote the set of permutations in $\S_n$ simultaneously avoiding
all patterns in $T$ by $\S_n(T)$.
\end{notation}

Permutations avoiding a 3-letter pattern were first considered in
\cite{Knuth}. In \cite{SS}, permutations and involutions avoiding
each set $T$ of 3-letter patterns were studied. Since then
restricted permutations and forbidden patterns were the subject of
many papers. One of the most frequently considered problems is the
enumeration of $\S_n(\tau)$ and $\S_n(T)$ for various patterns
$\tau$ and sets of patterns $T$. The inventory of cardinalities of
$|\S_n(T)|$ for $T\subseteq\S_3$ is given in \cite{SS}, and a
similar inventory for $|\S_n(\tau_1,\tau_2)|$, where
$\tau_1\in\S_3$ and $\tau_2\in\S_4$ is given in \cite{West2}. Some
results on $|\S_n(\tau_1,\tau_2)|$ for $\tau_1,\tau_2\in\S_4$ are
obtained in \cite{West1}. The exact formula for $|\S_n(1234)|$ and
the generating function for $|\S_n(12\dots k)|$ are found in
\cite{Gessel}. $|\S_n(1342)|=|\S_n(1423)|$ is obtained in
\cite{Bona}, and \cite{Stankova1,Stankova2} shows that
$|\S_n(3142)|=|\S_n(1342)|$. For a survey of results on pattern
avoidance, see \cite{KM}.

\begin{example} \label{ex:old} ~
\begin{itemize}
\item $|\S_n(\tau)|=C_n=\frac{1}{n+1}\binom{2n}{n}$, the $n$th
Catalan number, for any $\tau\in\S_3$.
\item $|\S_n(123,213)|=|\S_n(132,231)|=2^{n-1}$.
\item $|\S_n(123,132,213)|=F_n$, the $n$th Fibonacci number.
\item
$|\S_n(3142,2413)|=|\S_n(4132,4231)|=|\S_n(2431,4231)|=r_{n-1}$,
the $n$th large Schr\"oder number \cite[Sequence A006318]{Sloane},
given by $r_0=1$, $r_n=r_{n-1}+\sum_{j=0}^{n-1}{r_k r_{n-k}}$.
\end{itemize}
\end{example}

Another problem is finding (sets of) patterns $T_1$ and $T_2$ such
that $|\S_n(T_1)|=|\S_n(T_2)|$ for any $n\ge 0$. Such (sets of)
patterns are called \emph{Wilf-equivalent} and said to belong to
the same \emph{Wilf class}. There are three symmetry operations on
$\S_k$ that map every pattern onto a Wilf-equivalent
pattern:
\begin{itemize}

\item \emph{reversal} $r$: $r(\tau)(i)=\tau(n+1-i)$, i.e.
$r(\tau)$ is $\tau$ read right-to-left.

\item \emph{complement} $c$: $c(\tau)(i)=n+1-\tau(i)$, i.e.
$c(\tau)$ is $\tau$ read upside down.

\item $r\circ c=c\circ r$: $r\circ c(\tau)(i)=n+1-\tau(n+1-i)$,
i.e. $r\circ c(\tau)$ is $\tau$ read right-to-left upside down.
\end{itemize}

The set of patterns $\{\tau,r(\tau),c(\tau),r(c(\tau))\}$ is
called the \emph{symmetry class} of $\tau$.

\subsection{Dumont permutations} \label{subsec:dumont}

\begin{definition} \label{def:dumont1}
A \emph{Dumont permutation of the first kind} is a permutation
$\pi\in\S_{2n}$ of where each even entry is followed by a descent
and each odd entry is followed by an ascent or ends the string. In
other words, for every $i=1,2,\dots,2n$,
\[
\begin{split}
\pi(i) \text{ is even} &\implies i<2n \text{ and } \pi(i)>\pi(i+1),\\
\pi(i) \text{ is odd}  &\implies \pi(i)<\pi(i+1) \text{ or } i=2n.
\end{split}
\]
\end{definition}

\begin{definition} \label{def:dumont2}
A \emph{Dumont permutation of the second kind} is a permutation
$\pi\in\S_{2n}$ of where all entries at even positions are
deficiencies and all entries at odd positions are fixed points or
excedances. In other words, for every $i=1,2,\dots,n$,
\[
\begin{split}
\pi(2i)&<2i,\\
\pi(2i+1)&\ge \pi(2i+1).
\end{split}
\]
\end{definition}

\begin{notation} \label{not:dumont}
We denote the set of Dumont permutations of the first (resp.
second) kind of length $2n$ by $\D^1_{2n}$ (resp. $\D^2_{2n}$).
\end{notation}

\begin{example} \label{ex:dumont}
$\D^1_{2}=\D^2_{2}=\{21\}$, $\D^1_{4}=\{2143,3421,4213\}$,
$\D^2_{4}=\{2143,3142,4132\}$.
\end{example}

\begin{remark} \label{rem:odd}
Dumont permutations of odd length can be defined similarly to
those of even length. Then $\D^1_{2n+1}$ or $\D^2_{2n+1}$ are
obtained simply by adjoining $2n+1$ to the end of each permutation
in $\D^1_{2n}$ or $\D^2_{2n}$, respectively. Obviously,
$|\D^1_{2n+1}|=|\D^1_{2n}|$ and $|\D^2_{2n+1}|=|\D^2_{2n}|$.
\end{remark}

Dumont~\cite{Dumont} showed that
\[
|\D^1_{2n}|=|\D^2_{2n}|=G_{2n+2}=2(1-2^{2n+2})B_{2n+2},
\]
where $G_n$ is the $n$th Genocchi number, a multiple of the
Bernoulli number $B_n$. Lists of Dumont permutations $\D^1_{2n}$
and $\D^2_{2n}$ for $n\le 4$ as well as some basic information and
references for Genocchi numbers and Dumont permutations may be
obtained at \cite[A001469]{Sloane} and \cite{Ruskey}. We only note
that the exponential generating functions for the unsigned and
signed Genocchi numbers are given by
\[
\sum_{n=1}^{\infty}{G_{2n}\frac{x^{2n}}{(2n)!}}=x\tan\frac{x}{2},
\qquad \sum_{n=1}^{\infty}{(-1)^{n}G_{2n}\frac{x^{2n}}{(2n)!}}
=\frac{2x}{e^x+1}-x=-x\tanh\frac{x}{2}.
\]

Throughout this paper we will use the following obvious properties
of Dumont permutations.
\begin{fact} \label{fact:obvious1}
In any $\pi\in\D^1_{2n}$, $2$ is always followed by $1$, and
$2n-1$ is followed by $2n$ or is final in $\pi$.
\end{fact}
\begin{fact} \label{fact:obvious2}
In any $\pi\in\D^2_{2n}$, $\pi(2)=1$ and $\pi(2n-1)=2n-1 \text{ or
} 2n$.
\end{fact}

We define Wilf-equivalence and Wilf classes on Dumont permutations
in the same way as Wilf-equivalence on all permutations (and call
it $\D^1$-Wilf-equivalence or $\D^2$-Wilf-equivalence according to
the kind). Note that since reversal, complement, or reversal of
complement of Dumont permutations are no longer Dumont
permutations, it follows that permutations in the same symmetry
class are not necessarily Wilf-equivalent on Dumont permutations.

\begin{remark} \label{rem:dumont-like}
Sometimes slightly different permutations are defined as Dumont
permutations of either kind. Those will be useful later on, and we
will describe them now.

Permutations $\pi\in\S_{2n}$ in which each odd entry is followed
by ascent and each even entry is followed by a descent or ends the
string are obtained by applying the complement operation $c$ to
our Dumont permutations of the first kind, and we will call them
\emph{Dumont-like permutations of the first kind} and denote the
set of these in $\S_n$ by $\DD^1_n$.

Similarly, permutations $\pi\in\S_{2n}$ with $\pi(2i+1)>2i+1$ and
$\pi(2i)\le2i$ for all $i$ are obtained by applying the operation
$r\circ c=c\circ r$ to our Dumont permutations of the second kind,
and we will call them \emph{Dumont-like permutations of the second
kind} and denote the set of these in $\S_n$ by $\DD^2_n$.
\end{remark}

\subsection{Restricted Dumont permutations}

So far, there has been a single paper on restricted Dumont
permutations, namely Mansour~\cite{Mansour}. Most of it is devoted
to the study of 132-avoiding Dumont permutations of the first kind
(as there are no 132-avoiding Dumont permutations of the second
kind other than $21\in\D^2_2$). Specifically, 132-avoiding Dumont
permutations of the first kind which also avoid (contain exactly
once) certain other patterns $\tau$ are examined and enumerated.
However, Dumont permutations avoiding other patterns are briefly
examined as well.

\begin{notation} \label{not:dumont_avoid}
We denote the set of permutations in $\D^1_n$ or $\D^2_n$ avoiding
a pattern $\tau$ by $\D^1_n(\tau)$ or $\D^2_n(\tau)$,
respectively. If $T$ is a set of patterns, then we denote the set
of permutations in $\D^1_n$ or $\D^2_n$ avoiding all patterns in
$T$ by $\D^1_n(T)$ or $\D^2_n(T)$, respectively. We define
$\DD^1_n(T)$ and $\DD^2_n(T)$ (see Remark \ref{rem:dumont-like})
similarly.
\end{notation}

\begin{remark} \label{rem:symm}
Note that $\DD^1_{2n}(T)=c(\D^1_{2n}(c(T)))$ and
$\DD^2_{2n}(T)=r\circ c(\D^1_{2n}(r\circ c(T)))$, so
$|\DD^1_{2n}(T)|=|\D^1_{2n}(c(T))|$ and
$|\DD^2_{2n}(T)|=|\D^2_{2n}(r\circ c(T))|$.
\end{remark}

\begin{theorem} \label{thm:mansour}
{\rm (\cite[Theorems 2.2, 4.3]{Mansour})}
$|\D^1_{2n}(132)|=|\D^1_{2n}(231)|=|\D^1_{2n}(312)|=|\D^2_{2n}(321)|=C_n$.
Similarly,
$|\D^1_{2n+1}(132)|=|\D^1_{2n+1}(231)|=|\D^1_{2n+1}(312)|=|\D^2_{2n+1}(321)|=C_n$.
\end{theorem}

Another bit of notation will be useful before we proceed.
\begin{notation} \label{not:greater}
Let $\pi'$ and $\pi''$ be subsequences of a permutation $\pi$. We
say that $\pi'>\pi''$ if every entry of $\pi'$ is greater than
every entry of $\pi''$. Also, for a permutation $\pi$ and an
integer $m$, the string $\pi+m=m+\pi$ is obtained by adding $m$ to
every entry of $\pi$. We define $m-\pi$ similarly.
\end{notation}

All cardinalities may be obtained similarly by examining the
recursive structure of restricted permutations. For example, it
immediately follows from \cite[Proposition 2.1]{Mansour} that a
permutation in $\D^1_{2n}(132)$ can be of two types:
\begin{enumerate}
\item $\pi=(\pi',2n-1,2n,\pi'')$, where $\pi'>\pi''$ and
$\pi''\in\D^1_{2k}(132)$ and $\pi'-2k\in\D^1_{2n-2k-2}(132)$ for
some $1\le k\le n-1$.
\item $\pi=(2n,\pi',2n-1)$, where $\pi'\in\D^1_{2n-2}(132)$.
\end{enumerate}
(All permutations are given as lists in the one-line notation
unless otherwise indicated.)

The fact that $\pi'>\pi''$ follows from the fact that $\pi$ avoids
132. To see that $\pi''\in\D^1_{2k}$ (and not $\D^1_{2k+1}$) in
the first case, note that the minimum element of $\pi'$ must be
odd since it must be followed by an ascent.

Other patterns mentioned in Theorem \ref{thm:mansour} are treated
similarly.

In this paper, we use the same approach to analyze and enumerate
Dumont permutations of either kind avoiding certain patterns or
pairs of patterns of length 3 or 4. We will show that symmetry
operations on Dumont permutations do not necessarily produce a set
of patterns in the same Wilf class, hence we will only consider
certain cases in detail. We will see that the some cardinalities
of sets avoiding a given set of patterns are given by Catalan
numbers, powers of 2, little Schr\"oder numbers $s_n=r_n/2$, and
other known sequences.

We will frequently make use of Remark \ref{rem:symm}, especially
when considering $\D^1_{2n}(T)$ for a set of patterns $T=c(T)$,
such as $\{3142,2413\}$, $\{1342,4213\}$ or $\{1423,4132\}$, and
$\D^2_{2n}(T)$ for a pattern or a set of patterns $T=r\circ c(T)$,
such as $T=\{3142\}$.

\section{Dumont permutations avoiding 3-letter patterns}
\label{sec:3-letter}

As we mentioned before, Theorem \ref{thm:mansour} gives
$|\D^1_{2n}(132)|=|\D^1_{2n}(231)|=|\D^1_{2n}(312)|=|\D^2_{2n}(321)|=C_n$.

\begin{theorem} \label{thm:d1-213}
$|\D^1_{2n}(213)|=C_{n-1}$ for $n\ge 1$.
\end{theorem}

Note that $213=c(231)$ but $|\D^1_{2n}(213)|\ne|\D^1_{2n}(231)|$
whereas $132=c(312)$ and $|\D^1_{2n}(132)|=|\D^1_{2n}(312)|$.

\smallskip

\begin{proof}
From the observation in Fact \ref{fact:obvious1} it follows that a
permutation $\pi\in\D^1_{2n}(213)$ must end on $21$. Clearly,
$\D^1_2(213)=\{21\}$ (and hence $|\D^1_2(213)|=1=C_0$), so
consider $n\ge 2$. Let $\pi(1)=j$. Then $\pi=(j,\pi_1,\pi_2)$ for
some permutations $\pi_1>j>\pi_2$. Then $\pi_2=(\pi'',2,1)$, so
$j\ge 3$. Since $j$ is the minimal entry of $(j,\pi_1)$, it
follows that $j$ must be odd, i.e. $j=2k+1$ for some $k\ge 1$.
Since $(2k+1,\pi_1)>\pi_2$ and $(2k+1,\pi_1)$ starts with $2k+1$,
it follows that $(2k+1,\pi_1)$ ends with $2k+2$. Let
$(2k+1,\pi_1)=(2k+1,\pi',2k+2)$, then the last letter of $\pi'$ is
even, so $\pi'-(2k+2)\in\DD^1_{2n-2k-2}(213)$, in other words,
$c(\pi'-(2k+2))=2n+1-\pi'\in\D^1_{2n-2k-2}(231)$. Similarly,
$\pi_2=(\pi'',2,1)$, so the last letter of $\pi''$ is even and
hence $\pi''-2\in\DD^1_{2k-2}(213)$, i.e.
$c(\pi''-2)=2k+1-\pi''\in\D^1_{2k-2}(231)$.

Thus, the set $\D^1_{2n}(213)$ for $n\ge 1$ consists of all
permutations \[\pi=(2k+1,c(\rho')+2k+2,2k+2,c(\rho'')+2,2,1)\] for
some $k=1,2,\dots,n-1$, and some $\rho'\in\D^1_{2n-2k-2}(231)$ and
$\rho'\in\D^1_{2k-2}(231)$. Therefore, Theorem \ref{thm:mansour}
implies that, for $n\ge 2$,
\[
|\D^1_{2n}(213)|=\sum_{k=1}^{n-1}{C_{n-1-k}C_{k-1}}=C_{n-1}.
\]
Note that we could have done without the result of Theorem
\ref{thm:mansour} and proved our result using the recurrence
relation alone, but using Theorem \ref{thm:mansour} makes for a
slight shortcut in our argument.
\end{proof}

\begin{theorem} \label{thm:d2-231}
$|\D^2_{2n}(231)|=2^{n-1}$ for $n\ge 1$.
\end{theorem}

\begin{proof}
Since $\pi$ avoids 231, if follows that for any entry $j$ of
$\pi$, every element $<j$ to the left of $j$ must be lesser than
every element $<j$ to the right of $j$. From the observation in
Fact \ref{fact:obvious2} we know that if $\pi\in\D^2_{2n}(231)$,
then $\pi(2n-1)=2n$ or $\pi(2n-1)=2n-1$. We also have
$\pi(2n)<2n$. Therefore, $\pi(2n-1)=2n$ implies $\pi(2n)=2n-1$,
and $\pi(2n-1)=2n-1$ implies $\pi(2n)=2n-2$. In these two cases,
the last two entries of $\pi$ cannot be part of any occurrence of
231. Thus, $\pi\in\D^2_{2n}(231)$ if and only if either of the two
cases hold:
\begin{itemize}
    \item $\pi=(\pi',2n,2n-1)$ for any $\pi'\in\D^2_{2n-2}(231)$
    \item $\pi=(\widehat{\pi}',2n-1,2n-2)$ for any
    $\pi'\in\D^2_{2n-2}(231)$, where $\widehat{\pi}'$ obtains by
    replacing $2n-2$ with $2n$ in $\pi'$.
\end{itemize}
Therefore, $|\D^2_{2n}(231)|=2\cdot|\D^2_{2n-2}(231)|$ for $n\ge
2$, and $|\D^2_{2}(231)|=1$, so $|\D^2_{2n}(231)|=2^{n-1}$ for
$n\ge 1$.
\end{proof}

\smallskip

This recursive description allows us to determine the cycle
structure of permutations in $\D^2_{2n}(231)$. It is easy to see
inductively that any $\pi\in\D^2_{2n}(231)$ has $n$ cycles, each
of which contains exactly one odd entry and is of the form
$(2k-1)$ or $(2l,2l-2,\dots,2k,2k-1)$ for some $0\le k\le l\le n$.
Clearly, the above cycle must be followed by $l-k$ fixed points
$(2i-1)$, $k+1\le i\le l$. Thus, there is a natural bijection
between permutations $\D^2_{2n}(231)$ with $n-k$ fixed points and
weak $k$-compositions of $n$ (of which there are
$\binom{n-1}{k-1}$), where each cycle $(2l,2l-2,\dots,2k,2k-1)$
with $k\le l$ is mapped onto a part of size $l-k+1$ (and fixed
points are ``forgotten'').
\[
\begin{split}
\D^2_{8}(231)&\ni 21835476=(21)(8643)(5)(7) \mapsto 4=1+3\\
4=1+3 \mapsto &\underbrace{1}_{1}+\underbrace{3+0+0}_{3} \mapsto
(21)(8643)(5)(7)\in \D^2_{8}(231).
\end{split}
\]

\begin{theorem} \label{thm:d1-321}
$|\D^1_{2n}(321)|=1$ for $n\ge 0$.
\end{theorem}

\begin{proof}
We shall see that $\D^1_{2n}(321)=\{(2,1,4,3,\dots,2n,2n-1)\}$.
Indeed, if $\pi\in\D^1_{2n}(321)$ and $2=\pi(i)$ for some $i>1$,
then $\pi(i+1)=1$, so $\pi(i-1)>2$. Then the subsequence
$(\pi(i-1),\pi(i),\pi(i+1))=(\pi(i-1),2,1)$ of $\pi$ is an
occurrence of pattern 321. Hence, $\pi(1)=2$, so
$\pi=(2,1,\pi'+2)$ for some $\pi'\in\D^1_{2n-2}(321)$, and the
theorem follows by induction.
\end{proof}

\begin{theorem} \label{thm:d2-312}
$|\D^2_{2n}(312)|=1$ for $n\ge 0$.
\end{theorem}

\begin{proof}
We shall see that $\D^2_{2n}(312)=\{(2,1,4,3,\dots,2n,2n-1)\}$.
Indeed, if $\pi\in\D^2_{2n}(312)$ and $\pi(1)>2$, then $2=\pi(2i)$
for some $i\ge 2$, since $\pi(2)=1$. Then the subsequence
$(\pi(1),\pi(2),\pi(2i))=(\pi(1),1,2)$ of $\pi$ is an occurrence
of pattern 312. Therefore, $\pi(1)=2$, so $\pi=(2,1,\pi'+2)$ for
some $\pi'\in\D^2_{2n-2}(312)$, and the theorem follows by
induction.
\end{proof}

\begin{theorem} \label{thm:d1-123}
$|\D^1_{2n}(123)|=4$ for $n\ge 3$.
\end{theorem}

\begin{proof}
Note that $\D^1_{6}(123)=\{436215,562143,563421,564213\}$. If
$n\ge 4$ and $\pi\in\D^1_{2n}(123)$, then the subsequence of $\pi$
on letters $1,2,\dots,6$ belongs to $\D^1_{6}(123)$. Therefore,
$\pi$ cannot end on $2n-1$, so $2n-1$ is followed by $2n$ in
$\pi$, and hence $\pi=(2n-1,2n,\pi')$ for any
$\pi'\in\D^1_{2n-2}(123)$. This proves the theorem.
\end{proof}

\smallskip

Similarly, we can prove that
$\D^2_{2n}(123)=\D^2_{2n}(213)=\emptyset$ for $n\ge 3$, and
$\D^2_{2n}(132)=\emptyset$ for $n\ge 2$.

\smallskip

There is at most one Dumont permutation avoiding a given pair of
3-letter patterns simultaneously (except
$|\D^1_{2n}(123,132)|=2$). The following results are proved the
same way as Theorems \ref{thm:d1-321} and \ref{thm:d2-312}.

\begin{theorem} \label{thm:d1-d2-3-pairs}~\\[-26pt]
\[
\begin{split}
\D^1_{2n}(132,231)&=\{(2n,2n-2,\dots,4,2,1,3,\dots,2n-3,2n-1)\}\\
\D^1_{2n}(132,312)&=\{(\dots,2n-3,2n-2,4,3,2n-1,2n,2,1)\}\\
\D^1_{2n}(213,312)&=\{(3,5,\dots,2n-1,2n,\dots,6,4,2,1)\}\\
\D^1_{2n}(123,213)=\D^1_{2n}(132,213)&=\{(2n-1,2n,\dots,5,6,3,4,2,1)\}\\
\D^2_{2n}(231,321)=\D^1_{2n}(231,312)&=\{(2,1,4,3,\dots,2n,2n-1)\}\\
\D^1_{2n}(213,231)&=\emptyset \quad \text{for} \quad n\ge 2.
\end{split}
\]
\end{theorem}

\section{Dumont permutations avoiding 4-letter patterns}
\label{sec:4-letter}

When considering sets $\S_n(T)$ of permutations avoiding a given
set of patterns $T$, the first nontrivial cases arise when $T$ is
a single 3-letter permutation \cite{Knuth,SS}. In our case, the
first nontrivial restrictions by Dumont permutations are by
singletons in $\D^1_{4}=\{2143,3421,4213\}$ and
$\D^2_{4}=\{2143,3142,4132\}$. We were able to solve one of these
six cases.

\subsection{Avoiding a single pattern}\label{subsec:4-letter-single}

\begin{theorem} \label{thm:d2-3142}
$|\D^2_{2n}(3142)|=C_n$ for $n\ge 0$.
\end{theorem}

\begin{proof}
Let $\pi\in\D^2_{2n}(3142)$. Then $\pi(2)=1$. Consider two cases
based on parity of $\pi(1)$. Suppose that $\pi(1)=2k-1$ for some
$1\le k\le n$, then $\pi(2k-1)>2k-1$. If $\pi(i)<2k-1$ for some
$i>2k-1$, then $(\pi(1),\pi(2),\pi(2k-1),\pi(i))$ is an occurrence
of pattern 3142 in $\pi$. Therefore, $2\le \pi(i)\le 2k-2$ only if
$3\le i\le 2k-2$, i.e. all entries in $\{2,3,4\dots,2k-2\}$ must
occupy positions in $\{3,4,\dots,2k-2\}$, which is impossible.

Hence, we must have $\pi(1)=2k$ for some $1\le k\le n$. Since
$\pi(2k+1)>2k$, it follows, as before, that $2\le \pi(i)\le 2k-1$
only if $3\le i\le 2k$, i.e. all entries in $\{2,3,\dots,2k-1\}$
must occupy positions in $\{3,4,\dots,2k\}$. In other words,
$\pi=(2k,1,\pi_1+1,\pi_2+2k)$, where $\pi_1$ is a certain
permutation of $[2k-2]=\{1,2,\dots,2k-2\}$ avoiding pattern 3142,
and $\pi_2\in\D^2_{2n-2k}(3142)$. Furthermore, since $\pi_1+1$ is
a segment of $\pi\in\D^2_{2n}(3142)$ starting at position 3, it is
easy to see that $\pi_1\in\DD^2_{2k-2}(3142)$, so $\pi'=(r\circ c)
(\pi_1)\in\D^2_{2k-2}(3142)$ since $(r\circ c)(3142)=3142$. Thus,
$\pi\in\D^2_{2n}(3142)$ if and only if $\pi=(2k,1,(r\circ c)
(\pi')+1,\pi''+2k)$ for some $k=1,2,\dots,n$, and any
$\pi'\in\D^2_{2k-2}(3142)$ and $\pi''\in\D^2_{2n-2k}(3142)$. If we
let $a_n=|\D^2_{2n}(3142)|$, then $a_0=1$ and
$a_n=\sum_{k=1}^{n}{a_{k-1}a_{n-k}}$ for $n\ge 1$, so $a_n=C_n$.
\end{proof}

\smallskip

In fact, from the above proof, it is easy to see that any
$\pi\in\D^2_{2n}(3142)$ has the form
\[
\pi=(\pi_1,\pi_2+|\pi_1|,\pi_3+|\pi_2|+|\pi_1|,\dots,\pi_m+|\pi_{m-1}|+\dots+|\pi_1|)
\]
for some weak $m$-composition $n=n_1+n_2+\dots+n_m$ (i.e. all
$n_i>0$), and $\pi_i\in\D_{2n_i}(3142) \ (i=1,2,\dots,m)$ is of
the form
\[
\pi_i=(2n_i,1,(r\circ c)(\pi'_i)), \quad
\pi'_i\in\D_{2n_i-2}(3142).
\]
Furthermore, $(2n_1,2(n_1+n_2),2(n_1+n_2+n_3),\dots,2n)$ is the
sequence of left-to-right maxima of $\pi$. A natural bijective map
from $\D_{2n}(3142)$ to Dyck paths of $2n$ steps can be derived
recursively from the above decomposition, where left-to-right
maxima of $\pi$ (see above) are mapped to the steps leaving the
$x$-axis, and entries immediately following them are mapped to the
steps returning to the $x$-axis.

\bigskip

Since $321$ is a subpattern of $4132$, any $321$-avoiding
permutation also avoids $4132$. It is easy to see that
$\D^2_{2n}(321)$ consists of all permutations in $\D^2_{2n}$ in
which both the subsequence of entries in even positions and the
subsequence of entries in odd positions are increasing. It seems,
although we are currently unable to prove it\footnote{This result
has now been proved in \cite{BM-chebdumont}.}, that
$\D^2_{2n}(4132)$ also consists of these permutations only; in
other words, $|\D^2_{2n}(4132)|=|\D^2_{2n}(321)|=C_n$.

\begin{conjecture} \label{thm:d2-4132}
$|\D^2_{2n}(4132)|=C_n$ for $n\ge 0$.
\end{conjecture}

So far, we were unable to determine other cardinalities with
single pattern restrictions. The following lemma gives a
relationship between two such patterns, $4213$ and $1342=c(4213)$.

\begin{lemma} \label{lemma:1342<-4213}
If $A(x)$ and $B(x)$ are ordinary generating functions for
$|\D^1_{2n}(4213)|$ and $|\D^1_{2n}(1342)|$, respectively, then
\[
A(x)=\frac{1}{1-xB(x)}.
\]
\end{lemma}

\begin{proof}
Let $\pi\in\D^1_{2n}(4213)$. If $n>0$, then
$\pi=(\pi_1,2,1,\pi_2)$ for some $\pi_1<\pi_2$. Moreover, it easy
to see that $|\pi_1|+|\pi_2|=2n-2$, $\pi_1-2\in\DD^1(4213)$ so
$\pi'=|\pi_1|+3-\pi_1\in\D^1(1342)$, and
$\pi''=\pi_2-|\pi_1|\in\D^1(4213)$. It follows that
\[
A(x)=1+B(x)xA(x),
\]
which implies the lemma.
\end{proof}

Among single 4-letter Dumont permutations of either kind, only
$3142$ and $4132$ appear to be $\D^2$-Wilf-equivalent. It is easy
to see that $2143$ is not $\D^2$-Wilf-equivalent to 3142, and as
the following table shows, no two of the patterns in $\D^1_4$ are
$\D^1$-Wilf-equivalent.

\begin{center}
\begin{tabular}{c|cccccc}
  $n$ & 0 & 1 & 2 & 3 & 4 & 5 \\ \hline
  $|\D^1_{2n}(3421)|$ & 1 & 1 & 2 & 7 & 36 & 241 \\
  $|\D^1_{2n}(2143)|$ & 1 & 1 & 2 & 7 & 36 & 239 \\
  $|\D^1_{2n}(4213)|$ & 1 & 1 & 2 & 6 & 25 & 135 \\
\end{tabular}
\end{center}

\subsection{Avoiding a pair of patterns}\label{subsec:4-letter-pair}

\begin{theorem} \label{thm:d1-1342-1423}
$|\D^1_{2n}(1342,1423)|=s_{n+1}$ for $n\ge 0$.
\end{theorem}

Here $s_n$ is the $n$th little Schr\"oder number
\cite[A001003]{Sloane}, given by
$s_{n+1}=-s_{n}+2\sum_{k=1}^{n}{s_{k}s_{n-k}}$ ($n\ge 2$),
$s_1=1$, and the generating function
$s(x)=\sum_{n=1}^{\infty}{s_nx^n}=\left(1+x-\sqrt{1-6x+x^2}\right)/4$.

\smallskip

\begin{proof}
If $\pi\in\D^1_{2n}(1342,1423)$, then $\pi$ can be of two types:
\begin{itemize}
\item $\pi=(\pi_1,2n-1,2n,\pi_2)$ with $\pi_1=\emptyset$ or
$\pi_1>\pi_2\ne\emptyset$ (otherwise $\pi$ contains 1342); or,

\item $\pi=(\pi_1,2n,\pi_2,2n-1)$ with $\pi_1>\pi_2$ or
$\pi_1=\emptyset$ or $\pi_2=\emptyset$ (otherwise $\pi$ contains
1423).
\end{itemize}
In either case, if $\pi_1\ne\emptyset$, then the minimum entry of
$\pi_1$ is odd. Hence, the maximum entry of $\pi_2$ is even,
whether or not $\pi_1$ is nonempty. For $\pi$ of either type, the
entries $2n-1$ and $2n$ cannot be part of any occurrence of $1342$
or $1423$, so $\pi$ avoids 1342 and 1423 if and only if both
$\pi_1$ and $\pi_2$ avoid 1342 and 1423. Thus,
$\pi\in\D^1_{2n}(1342,1423)$ if and only if
\begin{itemize}
\item $\pi=(\pi'+2k,2n-1,2n,\pi'')$ for $1\le k\le n-1$,
$\pi'\in\D^1_{2n-2k-2}(1342,1423)$,\\
$\pi''\in\D^1_{2k}(1342,1423)$; or,

\item $\pi=(\pi'+2k,2n,\pi'',2n-1)$ for $0\le k\le n-1$,
$\pi'\in\D^1_{2n-2k-2}(1342,1423)$,\\
$\pi''\in\D^1_{2k}(1342,1423)$.
\end{itemize}
Therefore, for $a_n=|\D^1_{2n}(1342,1423)|$, we have $a_0=1$ and
$a_n=-a_{n-1}+2\sum_{k=0}^{n-1}{a_{k}a_{n-1-k}}$, so $\{a_n\}$
satisfies the same recurrence relation as $\{s_{n+1}\}$, and
hence, $a_n=s_{n+1}$.
\end{proof}

\smallskip

\begin{theorem} \label{thm:d1-2341-2413}
$|\D^1_{2n}(2341,2413)|=s_{n+1}$ for $n\ge 0$.
\end{theorem}

\begin{proof}
The proof is very similar to that of Theorem
\ref{thm:d1-1342-1423}. We get that $\pi\in\D^1_{2n}(2341,2413)$
if and only if
\begin{itemize}
\item $\pi=(\pi',2n-1,2n,\pi''+2k)$ for $0\le k\le n-2$,
$\pi'\in\D^1_{2k}(2341,2413)$,\\
$\pi''\in\D^1_{2n-2k-2}(2341,2413)$; or,

\item $\pi=(\pi',2n,\pi''+2k,2n-1)$ for $0\le k\le n-1$,
$\pi'\in\D^1_{2k}(2341,2413)$,\\
$\pi''\in\D^1_{2n-2k-2}(2341,2413)$,
\end{itemize}
Hence, the same recursive relation obtains, and the theorem
follows.
\end{proof}

\begin{theorem} \label{thm:d1-1342-2413}
$|\D^1_{2n}(1342,2413)|=s_{n+1}$ for $n\ge 0$.
\end{theorem}

\begin{proof}
Again, the proof is very similar to that of Theorem
\ref{thm:d1-1342-1423}. We have $\pi\in\D^1_{2n}(1342,2413)$ if
and only if
\begin{itemize}
\item $\pi=(\pi'+2k,2n-1,2n,\pi'')$ for $1\le k\le n-1$,
$\pi'\in\D^1_{2n-2k-2}(1342,2413)$,
$\pi''\in\D^1_{2k}(1342,2413)$; or,

\item $\pi=(\pi',2n,\pi''+2k,2n-1)$ for $0\le k\le n-1$,
$\pi'\in\D^1_{2k}(1342,2413)$,
$\pi''\in\D^1_{2n-2k-2}(1342,2413)$.
\end{itemize}
Hence, the same recursive relation obtains, and the theorem
follows.
\end{proof}

\begin{theorem} \label{thm:d1-2341-1423}
$|\D^1_{2n}(2341,1423)|=\frac{1}{\sqrt{17}}
\left(\left(\frac{3+\sqrt{17}}{2}\right)^n-\left(\frac{3-\sqrt{17}}{2}\right)^n\right)
=\left[\frac{1}{\sqrt{17}}
\left(\frac{3+\sqrt{17}}{2}\right)^n\right]$ for $n\ge 1$, where
$[a]$ is the integer closest to $a$.
\end{theorem}
Note that this sequence is \cite[A007482]{Sloane}, shifted right
by one position. In other words, $a_n=|\D^1_{2n}(2341,1423)|$ is
the number of subsets of $[2n-2]=\{1,2,\dots,2n-2\}$ where each
odd element $m$ has an even neighbor ($m-1$ or $m+1$). We also
have $a_0=1$, $a_1=1$, $a_2=3$, and $a_n=3a_{n-1}+2a_{n-2}$ for
$n\ge 3$.

\smallskip

\begin{proof}
As before, we have $\pi\in\D^1_{2n}(2341,1423)$ only if
\begin{enumerate}
\item $\pi=(\pi',2n-1,2n,\pi''+2k)$ for $0\le k\le n-2$,
$\pi'\in\D^1_{2k}(2341,1423)$,\\
$\pi''\in\D^1_{2n-2k-2}(2341,1423)$; or,
\item $\pi=(\pi'+2k,2n,\pi'',2n-1)$ for $0\le k\le n-1$,
$\pi'\in\D^1_{2n-2k-2}(2341,1423)$,\\
$\pi''\in\D^1_{2k}(2341,1423)$.
\end{enumerate}
However, now the entries $2n-1$ and $2n$ may be part of an
occurrence of a pattern $1423$ (in case 1) or a pattern $2341$ (in
case 2). In fact, it is easy to see that for $n\ge 3$ we must have
$k\in\{0,n-2\}$ in case 1, and $k\in\{0,n-2,n-1\}$ in case 2. In
other words, $\pi\in\D^1_{2n}(2341,1423)$ if and only if one of
the following holds:
\begin{enumerate}
\item $\pi=(2n-1,2n,\pi')$ for any
$\pi'\in\D^1_{2n-2}(2341,1423)$.

\item $\pi=(\pi',2n-1,2n,2n-2,2n-3)$ for any
$\pi'\in\D^1_{2n-4}(2341,1423)$.

\item $\pi=(2n,\pi',2n-1)$ for any
$\pi'\in\D^1_{2n-2}(2341,1423)$.

\item $\pi=(\pi',2n,2n-1)$ for any
$\pi'\in\D^1_{2n-2}(2341,1423)$.

\item $\pi=(2n-2,2n-3,2n,\pi',2n-1)$ for any
$\pi'\in\D^1_{2n-4}(2341,1423)$.
\end{enumerate}

Therefore, $a_n=3a_{n-1}+2a_{n-2}$ for $n\ge 3$, and the theorem
follows.
\end{proof}

\begin{theorem} \label{thm:d1-231-4213}
$\D^1_{2n}(231,4213)=\{(2,1,4,3,\dots,2n,2n-1)\}$ for $n\ge 1$.
\end{theorem}

\begin{proof}
For $n=1$, the result is obvious. Let $n\ge 2$ and
$\pi\in\D^1_{2n}(231,4213)$. If $\pi(2n)\ne 2n-1$, then $\pi$
contains a segment $(2n-1,2n,j)$ for some $j<2n-1$, which is an
instance of pattern $231$. Hence, $\pi(2n)=2n-1$. Now suppose that
$\pi(2n-1)\ne 2n$. If $2n-2$ occurs between $2n$ and $2n-1$, then
$2n-2$ must be followed by some $j<2n-2$, so $(2n,2n-2,j,2n-1)$ is
an instance of pattern $4213$. Therefore, $2n-2$ occurs before
$2n$, so $2n$ must be followed by some $j<2n-2$, and hence
$(2n-2,2n,j)$ is an instance of pattern $231$. Thus,
$\pi(2n-1)=2n$, in other words, $\pi=(\pi',2n,2n-1)$ for some
$\pi'\in\D^1_{2n-2}(231,4213)$. The rest is obvious.
\end{proof}

\begin{theorem} \label{thm:d1-1342-4213}
$|\D^1_{2n}(1342,4213)|=2^{n-1}$ for $n\ge 1$.
\end{theorem}

\begin{proof}
Let $\pi\in\D^1_{2n}(1342,4213)$, then $\pi=(\pi_1,2,1,\pi_2)$, so
$\pi_1<\pi_2$ if both $\pi_1,\pi_2\ne\emptyset$. Hence, the
largest letter in $\pi_1$ is even (since it must be followed by a
descent) and the last letter of $\pi_1$ is even (since it is
followed by 2). Also, $\pi_2$ must avoid $231$ since $(2,1,\pi_2)$
avoids $1342$. Therefore, $\pi=(\pi'+2,2,1,\pi''+2k)$ for some
$\pi'\in\DD^1_{2k-2}(1342,4213)$ (which means
$c(\pi')\in\D^1_{2k-2}(1342,4213)$) and
$\pi''\in\D^1_{2n-2k}(231,4213)$. Hence, either $\pi''=\emptyset$,
in which case $\pi=(2n+1-c(\bar\pi),2,1)$ for some
$\bar\pi\in\D^1_{2n-2}(1342,4213)$, or $\pi''\ne\emptyset$, in
which case $\pi=(\bar\pi,2n,2n-1)$ for some
$\bar\pi\in\D^1_{2n-2}(1342,4213)$. The rest is obvious.
\end{proof}

\begin{notation} \label{not:cx}
For the remaining part we will define
$C(x)=\frac{1-\sqrt{1-4x}}{2x}$, the ordinary generating function
for the sequence of Catalan numbers.
\end{notation}

\begin{theorem} \label{thm:d1-2413-3142}
The ordinary generating function for $|\D^1_{2n}(2413,3142)|$
$(n\ge 0)$ is given by
\[
F(x)=\frac{3-\sqrt{1-8x}}{2(1+x)}=\frac{1}{3}C\left(\frac{2(1+x)}{9}\right)
=\frac{1+2xC(2x)}{1+x}=\frac{1}{1-xC(2x)},
\]
so
\[
|\D^1_{2n}(2413,3142)|=(-1)^n+\sum_{k=1}^{n}{(-1)^{n-k}2^{k}C_{k-1}}
=\sum_{k=0}^{n-1}{\frac{n-k}{n+k}\binom{n+k}{n}2^k}+\delta_{n=0}
\]
is the convolution of ballot numbers and powers of 2.
\end{theorem}

Note also that $3142=c(2413)$, and that
$|\D^1_{2n}(2413,3142)|=C(2;n)$, the generalized Catalan number
\cite[A064062]{Sloane}.

\smallskip

\begin{proof}
Let $\pi\in\D^1_{2n}(2413,3142)$ ($n\ge 1$), and suppose that
$\pi(2n)=d=2k-1$ for some $1\le k\le n$. Consider a subsequence
$(a,b,c,d)$ of $\pi$. If $a<d$, $b>d$ and $c<d$, then $a<c$ since
$\pi$ avoids $2413$. Similarly, if $a>d$, $b<d$ and $c>d$, then
$a>c$ since $\pi$ avoids $3142$. Therefore,
\[
\pi=(\dots,\pi_6,\pi_5,\pi_4,\pi_3,\pi_2,2k-2,\pi_1,2k-1),
\]
where $2k-1>\pi_1>\pi_3>\pi_5>\dots$ and
$2k-2<\pi_2<\pi_4<\pi_6<\dots$. Also, $\pi_1$ and $\pi_2$ may be
empty, and while each $\pi_i$ for $i\ge 3$ must be nonempty, the
sequence $(\dots,\pi_6,\pi_5,\pi_4,\pi_3)$ may be empty. Note
that, if $i$ is odd (resp. even), then each nonempty $\pi_{i}$ is
followed by ascent (resp. descent), hence must end on an odd
(resp. even) number. Note also that the minima of all nonempty
$\pi_{2i-1}$ must be odd, while the maxima of all nonempty
$\pi_{2i}$ must be even. Therefore, $\pi_{2i-1}\in\D^1(2413,3142)$
for all $i\ge 1$, whereas $\pi_{2i}\in\DD^1(2413,3142)$, i.e.
$c(\pi_{2i})\in\D^1(2413,3142)$, for all $i\ge 1$. Finally, the
sum of the sizes of all $\pi_i$'s is $2n-2$.

Conversely, note that any permutation $\pi$ constructed as above
belongs to $\D^1_{2n}(2413,3142)$.

Let $a_n=|\D^1_{2n}(2413,3142)|$, and let
$F(x)=\sum_{n=0}^{\infty}{a_nx^n}$ be the ordinary generating
function for $\{a_n\}$, then the recursive structure of
permutations in $\D^1_{2n}(2413,3142)$ described above implies
that
\[
F(x)=1+xF(x)^2\frac{1}{1-(F(x)-1)},
\]
or, equivalently,
\[
(x+1)F(x)^2-3F(x)+2=0.
\]
This implies the theorem.
\end{proof}

\begin{theorem} \label{thm:d1-1423-4132}
The ordinary generating function for $|\D^1_{2n}(1423,4132)|$
($n\ge 0$) is given by
\[
G(x)=\frac{2-(1+x)C(x)}{2-x-(1+x)C(x)}=
\frac{1-3x-(1+x)\sqrt{1-4x}}{1-3x-2x^2-(1+x)\sqrt{1-4x}}.
\]
\end{theorem}

\begin{theorem} \label{thm:d1-2413-4132-1423-3142}
We have $|\D^1_{2n}(2413,4132)|=|\D^1_{2n}(1423,3142)|$ for $n\ge
0$ (i.e. $(2413,4132)$ and $(1423,3142)=(c(4132),c(2413))$ are
$\D^1$-Wilf-equivalent), and the ordinary generating function for
each sequence is given by
\[
H(x)=\frac{1+xC(x)-\sqrt{1-xC(x)-5x}}{2x(1+C(x))}.
\]
\end{theorem}

There are many directions in which to proceed further. We will
only mention several.

One such direction is to complete the investigation of single
forbidden patterns in $\D^1_4$ and $\D^2_4$, i.e. to find
$|\D^1_{2n}(\tau)|$ for $\tau=2143,3421,4213$ and
$|\D^2_{2n}(\tau)|$ for $\tau=2143,4132$. Another is to combine
the forbidden patterns of Section \ref{sec:4-letter} with
additional restrictions as in \cite{Mansour}. Yet another is to
find the complete distribution for the number of occurrences of
these patterns possibly combined with other restrictions, or to
find equidistributed statistics on some of these restricted sets.
Finally, the restrictions that define $\D^1$ may be generalized to
strings with repeated letters. It remains to be seen if a
generalization to words is possible for $\D^2$.

\end{document}